\documentstyle[12pt]{article}
\def\proof{\noindent{\bf Proof:}\hskip10pt}
\def\QED{\hfill\vrule height 1.5ex width 1.4ex depth -.1ex \vskip20pt}

\font\tenmath=msbm10 scaled 1200
\font\sevenmath=msbm7 scaled 1200
\font\fivemath=msbm5 scaled 1200
\newfam\mathfam \textfont\mathfam=\tenmath
\scriptfont\mathfam=\sevenmath \scriptscriptfont\mathfam=\fivemath
\def\math{\fam\mathfam}
\begin{document}
\def \\ { \cr }
\def\R{{\math R}}
\def\N{{\math N}}
\def\E{{\math E}}
\def\P{{\math P}}
\def\Z{{\math Z}}
\def\Q{{\math Q}}
\def\C{{\math C}}
\def \e{{\rm e}}
\def \f{{\cal F}}
\def \g{{\cal G}}
\def \h{{\cal H}}
\def \d{{\tt d}}
\def \k{{\tt k}}
\def \i{{\tt i}}
\def \B{{\bf B}}
\def \L{{\cal L}}
\newcommand{\ed}{\mbox{$ \ \stackrel{d}{=}$ }}
\newtheorem{theorem}{Theorem}
\newtheorem{proposition}[theorem]{Proposition}
\newtheorem{lemma}[theorem]{Lemma}
\newtheorem{corollary}[theorem]{Corollary}
\newtheorem{definition}[theorem]{Definition}
\centerline{\Large \bf {Some Connections Between 
(Sub)Critical}}
\vskip 2mm
\centerline{\Large \bf  Branching Mechanisms 
and Bernstein Functions}

\vskip 1cm
\centerline{\Large \bf Jean Bertoin$^{(1)}$, Bernard Roynette$^{(2)}$, and
Marc Yor$^{(1)}$}

\vskip 1cm
\noindent
\noindent
(1) {\sl Laboratoire de Probabilit\'es et Mod\` eles Al\'eatoires
and Institut universitaire de France,
 Universit\'e Pierre et Marie Curie,  175, rue du Chevaleret,
 F-75013 Paris, France.}
\vskip 2mm
\noindent
(2) {\sl {Institut Elie Cartan, Campus Scientifique,
BP 239, Vandoeuvre-l\` es-Nancy Cedex
F-54056, France
}}

\begin{abstract}
We describe some connections, via composition, between two functional
spaces: the space of (sub)critical branching mechanisms and the space of
Bernstein functions. The functions ${\bf e}_\alpha: x\to x^{\alpha}$ where
$x\geq0$ and
$0<\alpha\leq 1/2$, and in particular the critical parameter
$\alpha=1/2$, play a distinguished role.
\end{abstract}

\section{Introduction} This note is a prolongation of \cite{RY} where the
following remarkable property of the function ${\bf e}_\alpha: x\to
x^{\alpha}$  was pointed at for $\alpha=1/2$: if $\Psi$ is a (sub)critical
branching mechanism, then $\Psi\circ {\bf e}_{1/2} $ is a Bernstein
function (see the next section for the definition of these notions). In
the present work, we first show that this property extends to every
$\alpha\in]0,1/2]$. Then we characterize the class of so-called internal
functions, i.e. that of Bernstein functions $\Phi$ such that the compound
function
$\Psi\circ \Phi$ is again a Bernstein function for every 
(sub)critical branching mechanism $\Psi$. In the final section, 
we gather classical results on transformations of completely monotone
functions, Bernstein functions and (sub)critical branching mechanisms
which are used in our analysis.

\section{Some functional spaces}

\subsection{Completely monotone functions}

For every   Radon measure $\mu$ 
on $[0,\infty[$, we associate the function
$\L_{\mu}: ]0,\infty[\to[0,\infty]$ defined by
\begin{equation}\label{eq0}
\L_{\mu}(q):=\int_{]0,\infty[}\e^{-qx}\mu(dx)\,,
\end{equation}
i.e. $\L_{\mu}$ is the Laplace transform of $\mu$.
We denote by 
\begin{equation}\label{eq1}
{\bf CM}:=\left\{\L_\mu: \L_\mu(q)<\infty\hbox{ for all }q>0\right\}\,,
\end{equation}
which is an algebraic convex cone (i.e. a convex cone which is
further stable under inner product). The celebrated theorem of Bernstein
(see for instance Theorem 3.8.13 in \cite{Jacob}) identifies
${\bf CM}$ with the space of completely monotone functions, i.e.
functions $f: ]0,\infty[
\to[0,\infty[$ of class ${\cal C}^{\infty}$ such
that for every integer $n\geq1$, the $n$-th derivative $f^{(n)}$ of $f$
has the same sign as $(-1)^n$.
Recall from monotone convergence that $\L_\mu$ has a (possibly infinite) limit at $0+$
which coincides with the total mass of $\mu$.

We shall focus on two natural sub-cones of ${\bf CM}$: 
\begin{equation}\label{eq2}
\B_1:=\left\{\L_{\mu}: \int_{]0,\infty[}(1\wedge
x^{-1})\mu(dx)<\infty\right\}
\end{equation}
We further denote by $\B_1^{\downarrow}$ the sub-space of functions in
$\B_1$ which are the Laplace transforms of absolutely
continuous measures with a decreasing density :
\begin{equation}\label{eq4}
\B_1^{\downarrow}:=\left\{\L_\mu: \mu(dx)=g(x)dx, g \hbox{ decreasing 
and } 
\int_{0}^{\infty}(1\wedge x^{-1})g(x)dx<\infty\right\}.
\end{equation}
Note that the density $g$ then has limit $0$ at infinity.

\subsection{Bernstein functions}
For every triple $(a,b,\Lambda)$ with $a,b\geq0$ and $\Lambda$ a
positive measure on
$]0,\infty[$ such that
\begin{equation}\label{eq5}
\int_{]0,\infty[}(x\wedge 1) \Lambda(dx)<\infty\,,
\end{equation}
we associate the function $\Phi_{a,b,\Lambda}: ]0,\infty[\to [0,\infty[$
defined by
\begin{equation}\label{eq6}
\Phi_{a,b,\Lambda}(q):=a +bq+\int_{]0,\infty[}(
1-\e^{-q x})\Lambda(dx)\,,
\end{equation}
and call  $\Phi_{a,b,\Lambda}$ the Bernstein function with characteristics
$(a,b,\Lambda)$.
We denote the convex cone of Bernstein functions by
\begin{equation}\label{eq7}
\B_2:=\left\{\Phi_{a,b,\Lambda}: a,b\geq0\hbox{ and $\Lambda$
positive measure fulfilling (\ref{eq5})}\right\}.
\end{equation}

It is well-known that $\B_2$ can be identified
with the space of real-valued ${\cal C}^{\infty}$ functions $f:
]0,\infty[\to[0,\infty[$ such that for every integer $n\geq1$, the $n$-th
derivative $f^{(n)}$ of $f$ has the same sign as $(-1)^{n-1}$. See
Definition 3.9.1 and Theorem 3.9.4 in \cite{Jacob}.

Bernstein functions appear as Laplace exponents of subordinators,
 see e.g. Chapter 1 in \cite{Besf},
Chapter 6 in \cite{Sato}, or Section 3.9 in \cite{Jacob}. This means that
$\Phi\in\B_2$ if and only if there exists an increasing process
$\sigma=(\sigma_t, t\geq0)$ with values in $[0,\infty]$ ($\infty$ serves
as absorbing state) with independent and stationary increments as long as
$\sigma_t<\infty$, such that for every $t\geq0$
$$\E(\exp(-q\sigma_t))\,=\,\exp(-t\Phi(q))\,,\qquad q>0.$$
In this setting, $a$ is known as the killing rate, $b$ as the drift coefficient, and 
$\Lambda$ as the L\'evy measure.

We shall further denote by $\B_2^{\downarrow}$ the subspace of
Bernstein functions  for which
the L\'evy measure  is absolutely continuous with a monotone
decreasing density, viz.
$$\B_2^{\downarrow}:=\left\{\Phi_{a,b,\Lambda}: a,b\geq0
\hbox{ and }\Lambda(dx)=g(x)dx, g\geq0
\hbox{ decreasing  and } 
\int_{0}^{\infty}(x\wedge 1)g(x)dx<\infty\right\}.$$

\subsection{ (Sub)critical branching mechanisms}
For every triple $(a,b,\Pi)$ with $a,b\geq0$ and $\Pi$ positive measure on
$]0,\infty[$ such that
\begin{equation}\label{eq8}
\int_{]0,\infty[}(x\wedge x^2) \Pi(dx)<\infty
\end{equation}
we associate the function
$\Psi_{a,b,\Pi}: ]0,\infty[\to[0,\infty[$ defined by
\begin{equation}\label{eq9}
\Psi_{a,b,\Pi}(q):=a q +bq^2+\int_{]0,\infty[}(
\e^{-q x}-1 +q x)\Pi(dx)\,,
\end{equation}
and denote the convex cone of such functions by
\begin{equation}\label{eq10}
\B_3:=\left\{\Psi_{a,b,\Pi}: a,b\geq0\hbox{ and $\Pi$ a positive measure
such that (\ref{eq8}) holds}\right\}
\end{equation}

Functions in $\B_3$ are convex increasing functions of class ${\cal
C}^{\infty}$ that vanish at $0$; they coincide with the class
of branching mechanisms for (sub)critical continuous state branching
processes, where  (sub)critical means critical or sub-critical. See Le
Gall \cite{LG} on page 132.

 Alternatively, functions in the space $\B_3$ can also be viewed as 
Laplace exponents
 of L\'evy processes with no positive jumps that do not drift to
$-\infty$ (or, equivalently, with nonnegative mean). In this setting, $a$
is the drift coefficient,
$2b$ the Gaussian coefficient, and $\Pi$ the image of the L\'evy measure
by the map $x\to-x$.
 See e.g. Chapter VII in \cite{Belp}.

\section{Composition with ${\bf e}_{\alpha}$}
Stable subordinators correspond to a remarkable one-parameter family
of Bernstein functions denoted here by
$({\bf e}_{\alpha}, 0<\alpha<1)$, where
$$
{\bf e}_{\alpha}(q):=q^{\alpha}\,=\,
{\alpha\over
\Gamma(1-\alpha)}
\int_{0}^{\infty}
(1-\e^{-q x}) x^{-1-\alpha}dx\,, \qquad q>0\,.
$$

\begin{theorem}\label{T1} The following assertions are equivalent:

\noindent {\rm (i)} $\alpha\in]0,1/2]$.

\noindent {\rm (ii)} For every $\Psi\in \B_3$, $\Psi\circ {\bf
e}_{\alpha}\in\B_2$.
\end{theorem}

The implication (ii) $\Rightarrow$ (i) is immediate.
Indeed, $\Psi_{0,1,0}: q\to q^2$ belongs to $\B_3$, but
${\bf e}_{2\alpha}=\Psi_{0,1,0}\circ {\bf e}_{\alpha}$ is in
$\B_2$ if and only if $2\alpha\leq 1$.
However, the converse (i) $\Rightarrow$ (ii) is not straightforward and
relies on the following technical lemma, which appears as Lemma VI.1.2 in
\cite{RY}. Here, for the sake of completeness, we provide  a proof.
\begin{lemma}\label{L1} For $\alpha\in ]0,1/2]$,
let $\sigma^{(\alpha)}=(\sigma^{(\alpha)}_x, x\geq0)$ be a stable subordinator
with index $\alpha$ with Laplace transform 
$$\E\left(\exp\left(-{q} \sigma^{(\alpha)}_x\right)\right)
=\exp(-x q^{\alpha})\,,\qquad x,q>0\,.$$
Denote by $p^{(\alpha)}(x,t)$ the density of the law of 
$\sigma^{(\alpha)}_x$. Then for every $,x,t>0$, we have
$$p^{(\alpha)}(x,t)\leq  {\alpha\over \Gamma(1-\alpha)} x
t^{-(1+\alpha)}\,.$$
\end{lemma} 

\noindent{\bf Remark :} The bound in Lemma \ref{L1} is sharp, as it is
well-known that for any $0<\alpha<1$ and each fixed $t>0$ 
$$p^{(\alpha)}(x,t) \,\sim\,{\alpha\over \Gamma(1-\alpha)} x
t^{-(1+\alpha)}\,,\qquad x\to\infty.$$
More precisely, there is a series representation of $p^{(\alpha)}(x,t)$, see Formula
(2.4.7) on page 90
in Zolotarev \cite{Zol}:
 $$p^{(\alpha)}(x,1)={1\over \pi}\sum_{n=1}^{\infty}(-1)^{n-1}{\Gamma(n\alpha+1)\over \Gamma(n+1))}\sin(\pi n \alpha) x^{-n\alpha -1}\,.$$
 Using the identity 
 $$\Gamma(\alpha)\Gamma(1-\alpha)={\pi \over\sin (\alpha \pi)}\,,$$
this agrees  of course with the above estimate. It is interesting to note
that the second leading term in the expansion, 
$$-{\Gamma(2\alpha+1)\over 2\pi}\sin(2\pi \alpha) x^{-2\alpha -1},$$
is negative for $\alpha<1/2$, but positive for $\alpha >1/2$.
So the bound in Lemma \ref{L1} would fail for $\alpha>1/2$.

\proof In the case $\alpha=1/2$, there is an explicit expression for the
density
$$p^{(1/2)}(x,t)\,=\,{x\over 2\sqrt{\pi t^3}}\exp\left(-{x^2\over
4t}\right)\,,$$
from which the claim is obvious (recall that $\Gamma(1/2)=\sqrt \pi$).

In the case $\alpha<1/2$, we start from the identity
$$\exp(-xq^{\alpha})
\,=\,\int_{0}^{\infty}\e^{-qt} p^{(\alpha)}(x,t)dt
\,, $$
and take the derivative in the variable $q$ to get
$$\alpha q^{\alpha-1}\exp(-xq^{\alpha})
\,=\,\int_{0}^{\infty}\e^{-qt} {t\over x}
p^{(\alpha)}(x,t)dt\,,$$
and then
$$\alpha q^{\alpha-1}\left(1-\exp(-xq^{\alpha})\right)
\,=\,\int_{0}^{\infty}\e^{-qt}\left({\alpha\over \Gamma(1-\alpha)}
t^{-\alpha}- {t\over x} p^{(\alpha)}(x,t)\right)dt\,.$$

Denote the left hand-side by $g(x,q)$, and take the derivative in the
variable $x$. We obtain
$${\partial g(x,q)\over \partial x}
\,=\,\alpha q^{2\alpha-1}\e^{-xq^{\alpha}}
\,=\,\alpha q^{2\alpha-1}\int_{0}^{\infty}\e^{-qt}
p^{(\alpha)}(x,t)dt\,.$$
On the other hand, since $1-2\alpha>0$,
$$q^{2\alpha-1}\,=\,{1\over
\Gamma(1-2\alpha)}\int_{0}^{\infty}\e^{-qs}s^{-2\alpha}ds\,,$$
and hence
$${\partial g(x,q)\over \partial x}
\,=\,{\alpha\over
\Gamma(1-2\alpha)}\int_{0}^{\infty}{ds\over
s^{2\alpha}}\int_{0}^{\infty}dt
\e^{-q(s+t)}
 p^{(\alpha)}(x,t)\,.$$
The change of variables $u=t+s$ yields
$${\partial g(x,q)\over \partial x}
\,=\,{\alpha\over
\Gamma(1-2\alpha)}\int_{0}^{\infty}du\int_{0}^{u}{ds\over
s^{2\alpha}}
\e^{-qu} 
 p^{(\alpha)}(x,u-s)\,;$$
and since $g(0,t)=0$, we finally obtain the identity
\begin{eqnarray*}
& &\int_{0}^{\infty}\e^{-qt}\left({\alpha\over \Gamma(1-\alpha)}
t^{-\alpha}- {t\over x} p^{(\alpha)}(x,t)\right)dt\\
&=&
{\alpha\over
\Gamma(1-2\alpha)}\int_{0}^{x}dy\int_{0}^{\infty}du\int_{0}^{u}{ds\over
s^{2\alpha}}
\e^{-qu} 
 p^{(\alpha)}(x,u-s)\,.
\end{eqnarray*}
Inverting the Laplace transform, we conclude that
$${\alpha\over \Gamma(1-\alpha)}
t^{-\alpha}- {t\over x} p^{(\alpha)}(x,t)
\,=\,
{\alpha\over
\Gamma(1-2\alpha)}\int_{0}^{x}dy\int_{0}^{t}{ds\over
s^{2\alpha}}
 p^{(\alpha)}(x,t-s)\,,$$
which entails our claim. \QED

We are now able to prove Theorem \ref{T1}.

\proof Let $\Psi_{a,b,\Pi}\in\B_3$.
Since both $a{\bf e}_{\alpha}$ and $b{\bf e}_{2\alpha}$ are Bernstein
functions, there is no loss of generality in assuming that
$a=b=0$. Set for $t>0$
$$\nu_{\alpha}(t):=
{\alpha\over
\Gamma(1-\alpha) t^{1+\alpha}}\int_{0}^{\infty}\Pi(dx) x
\left(1-{\Gamma(1-\alpha) t^{1+\alpha}\over \alpha x}
 p^{(\alpha)}(x,t)\right)\,.$$
It follows from Lemma \ref{L1} that $\nu_{\alpha}(t)\geq0$.
We have for every $q>0$
\begin{eqnarray*}
& &\int_{0}^{\infty}(1-\e^{-qt})\nu_{\alpha}(t)dt\\
&=&\int_{0}^{\infty}\Pi(dx) x \int_{0}^{\infty}dt
\left({\alpha (1-\e^{-qt})\over
\Gamma(1-\alpha) t^{1+\alpha}}
-{p^{(\alpha)}(x,t)\over x}+
\e^{-qt}
{p^{(\alpha)}(x,t)\over x}\right)\\
&=&\,\int_{0}^{\infty}\Pi(dx) x\left(q^{\alpha}-{1\over
x}+{\e^{-q^{\alpha}x}\over x}\right)\\
&=& \Psi_{0,0,\Pi}({\bf e}_{\alpha}(q))\,.
\end{eqnarray*}
As this quantity is finite for every $q>0$, this shows that 
$\Psi_{0,0,\Pi}\circ{\bf e}_{\alpha}\in\B_2$. \QED

\noindent {\bf Remark :} The proof gives a stronger
result than that stated in Theorem \ref{T1}. Indeed, we specified the
L\'evy measure
$\nu_{\alpha}$ of $\Psi_{0,0,\Pi}\circ{\bf e}_{\alpha}$. 
Furthermore, in the case $\alpha=1/2$, this expression shows that
$\Psi_{0,0,\Pi}\circ{\bf
 e}_{1/2}\in\B_2^{\downarrow}$. It is interesting to combine this
observation with the forthcoming Proposition \ref{P2} : for every
$\Psi\in\B_3$, $\Psi\circ{\bf
 e}_{1/2}\in\B_2^{\downarrow}$, thus ${\rm Id}\times (\Psi\circ{\bf
 e}_{1/2}): q\to q\Phi(\sqrt q)$ is again in $\B_3$, and in turn  
${\bf e}_{1/2}\times (\Psi\circ{\bf
 e}_{1/4})\in\B_2^{\downarrow}$. More generally, we have by iteration
that for every integer $n$
$${\bf e}_{2-2^{1-n}}\times (\Psi\circ{\bf
 e}_{2^{-n}})\in\B_3\,,$$
and 
$${\bf e}_{1-2^{-n}}\times (\Psi\circ{\bf
 e}_{2^{-n-1}})\in\B_2^{\downarrow}\,.$$

\section{Internal functions}
It is well-known that the cone ${\bf CM}$ of completely monotone
functions and the cone $\B_2$ of Bernstein functions are both stable by
right composition with a Bernstein function; see Proposition \ref{P3}
below. Theorem
\ref{T1} incites us to consider also compositions
of (sub)critical branching mechanisms and Bernstein functions; we make
the following definition :

\begin{definition} A Bernstein function $\Phi\in\B_2$ is said
{\rm internal} if $\Psi\circ \Phi\in\B_2$ for every $\Psi\in\B_3$.
\end{definition}

Theorem \ref{T1} shows that the functions ${\bf e}_{\alpha}$ are internal
if and only if $\alpha\in]0,1/2]$. The critical parameter $\alpha=1/2$
plays a distinguished role. Indeed, we could also prove Theorem \ref{T1}
using the following alternative route. First, we check that ${\bf
e}_{1/2}$ is internal (see \cite{RY}), and then we deduce by subordination
that for every
$\alpha<1/2$ that $\Psi\circ {\bf e}_{\alpha}=
\Psi\circ {\bf e}_{1/2}\circ {\bf e}_{2\alpha}$ is again a Bernstein
function for every $\Psi\in\B_3$. Developing this argument, we easily
arrive at the following characterization of internal functions :

\begin{theorem}\label{T2} Let $\Phi=\Phi_{a,b,\Lambda}\in\B_2$
be a Bernstein function. The following assertions are then equivalent:

\noindent{ \rm (i)} $\Phi$ is internal,

\noindent{ \rm (ii)} $\Phi^2\in\B_2$,

\noindent{ \rm (iii)} $b=0$ and there exists a subordinator $\sigma=
(\sigma_t, t\geq0)$ such that
$$\Lambda(dx)\,=\,c\int_{0}^{\infty} t^{-3/2} \P(\sigma_t\in dx)dt\,.$$

\end{theorem}

\proof (i) $\Rightarrow$ (ii) is obvious as $\Psi_{0,1,0}\circ \Phi
=\Phi^2$.

(ii) $\Rightarrow$ (i). We know from Theorem \ref{T1} or \cite{RY} that
for every $\Psi\in\B_3$, $\Psi\circ {\bf e}_{1/2}\in\B_2$.
It follows by subordination that for every Bernstein function
$\kappa\in\B_2$, $\Psi\circ {\bf e}_{1/2}\circ \kappa\in\B_2$.
Take $\kappa=\Phi^2$, so $ {\bf e}_{1/2}\circ \kappa=\Phi$,
and hence $\Phi$ is internal.

(iii) $\Rightarrow$ (ii) Let $\kappa$ denote the Bernstein function of
$\sigma$. We have
\begin{eqnarray*}
\Phi(q)\,&=&\,a +\int_{]0,\infty[}(1-\e^{-qx})\Lambda(dx)\\
&=&\,a + c \int_{]0,\infty[}\int_{0}^{\infty}dt (1-\e^{-qx}) t^{-3/2}
\P(\sigma_t\in dx)\\
&=&\,a + c \int_{0}^{\infty}dt (1-\e^{-t\kappa(q)}) t^{-3/2}\,.
\end{eqnarray*}
The change of variables $t\kappa(q)=u$ yields
$$\Phi(q)\,=\,a+c'\sqrt{\kappa(q)}$$
and hence
$$\Phi^2(q)\,=\,a^2 + 2ac'\sqrt{\kappa(q)} +c'^2\kappa(q)\,.$$
Since $\kappa^{1/2}={\bf e}_{1/2}\circ \kappa$ is again a Bernstein
function, we thus see that $\Phi^2\in\B_2$.

(ii) $\Rightarrow$ (iii) Recall that the drift coefficient $b$ of
$\Phi_{a,b,\Lambda}$ is given by
$$\lim_{q\to\infty} \Phi_{a,b,\Lambda}(q)/q\,=\,b\,;$$
see e.g. page 7 in \cite{Besf}. It follows immediately that
$b=0$ whenever $\kappa:=\Phi_{a,b,\Lambda}^2\in\B_2$.
Recall from Sato \cite{Sato} on page 197-8 that if $\tau^{(1)}$ and
$\tau^{(2)}$ are two independent subordinators with respective Bernstein
functions $\Phi^{(1)}$ and $\Phi^{(2)}$, then the compound process
$\tau^{(1)}\circ \tau^{(2)}:=\tau^{(3)}$ is again a subordinator with
Bernstein function $\Phi^{(3)}:=\Phi^{(2)}\circ \Phi^{(1)}$; moreover
its L\'evy measure $\Lambda^{(3)}$ is given by
$$\Lambda^{(3)}(dx)=\int_{0}^{\infty}\P(\tau^{(1)}_t\in
dx)\Lambda^{(2)}(dt)\,,$$
where $\Lambda^{(2)}$ denotes the L\'evy measure of $\tau^{(2)}$.
As $\Phi_{a,b,\Lambda}={\bf e}_{1/2}\circ \kappa$, and the L\'evy measure
of ${\bf e}_{1/2}$ is $ct^{-3/2}dt$ with $c=1/(2\sqrt \pi)$, 
we deduce  that 
$$\Lambda(dx)\,=\,c\int_{0}^{\infty} \P(\sigma_t\in dx) t^{-3/2}dt\,.$$

The proof of Theorem \ref{T2} is now complete. \QED

It is noteworthy that if $\Phi_{a,b,\Lambda}$ is internal
and $\Lambda\not\equiv 0$,
then
$$\int_{]0,\infty}x \Lambda(dx)\,=\,\infty\,.$$
Indeed,
$$\int_{]0,\infty}x \Lambda(dx)
\,=\,c\int_{0}^{\infty}\int_{]0,\infty[} x \P(\sigma_t\in dx)
t^{-3/2}dt
\,=\,c\int_{0}^{\infty}\E(\sigma_1)
t^{-1/2}dt\,=\,\infty\,.$$
For instance, the Bernstein function $q\to \log(1+q)$ of the gamma
subordinator is not internal.

\begin{corollary}\label{C1} For every $\Psi\in\B_3$, 
wewrite $\Phi$ for the inverse function of $\Psi$
and then $\Phi'$ for its derivative. Then
$1/\Phi'$ is internal.
\end{corollary}

\proof
It is known (see Corollary \ref{C2} below) that
$1/\Phi'$ is a Bernstein function; let us check that its
square is also a Bernstein function.

We know that $\Psi''\in \B_1$ (Proposition \ref{P1} below)
 and
$\Phi\in \B_2$ (Proposition \ref{P4} below); we deduce
from Proposition \ref{P3} that $\Psi''\circ \Phi\in\B_1$.
If we write $I(f): x\to\int_{0}^{x}f(y)dy$ for
every locally integrable function $f$, then again by Proposition \ref{P1},
we get that  $I(\Psi''\circ \Phi)$
is a Bernstein function.

Now
$$\Psi''\,=\,-{\Phi''\circ \Psi\over
(\Phi'\circ \Psi)^3}\,,$$
so 
$$\Psi''\circ \Phi\,=\,-{\Phi''\over
(\Phi')^3}\,,$$
and we conclude that
$${1\over 2(\Phi')^2}\,=\,I(\Psi''\circ \Phi)\in
\B_2\,.$$
\QED

\section{Some classical results and their consequences}
For convenience, this section gathers some classical transformations
involving $\B_j$, $j\in\{1,2,3\}$ and related subspaces, which
have been used in the preceding section. We start by considering
derivatives and indefinite integrals. The following statement is
immediate.

\begin{proposition}\label{P1} Let $j=2,3$ and $f:]0,\infty[\to[0,\infty[$
be a ${\cal C}^{\infty}$-function with derivative
$f'$. For $j=3$, we further suppose that 
$\lim_{q\to0}f(q)=0$.
There is the equivalence

$$f\in\B_j \ \Longleftrightarrow \ f'\in \B_{j-1}\,.$$
\end{proposition}

The next statement is easily checked using integration by parts.

\begin{proposition}\label{P2}  Let $j=2,3$ and consider two functions
$f,g:]0,\infty[\to[0,\infty[$ which are related by the identity
$f(q)=qg(q)$. 
Then there is the equivalence 

$$f\in\B_j \hbox{ and } \lim_{q\to0} f(q)=0 \ \Longleftrightarrow \ g\in
\B_{j-1}^{\downarrow}
\,.$$
\end{proposition}
Proposition \ref{P2} has well-known probabilistic interpretations.
First, let $\sigma$ be a subordinator with Bernstein function $f\in \B_2$
with unit mean, viz. $\E(\sigma_1)=1$, which is equivalent to $f'(0+)=1$.
Then the completely monotone function $g(q):=f(q)/q$ 
 is the Laplace transform of a probability measure on $\R_+$. The
latter appears in the renewal theorem for subordinators (see e.g.
\cite{BvHS}); in particular it describes the weak limit of the so-called
age process
$A(t)=t-g_t$ as
$t\to\infty$, where
$g_t:=\sup\left\{\sigma_s: \sigma_s<t\right\}.$
Second, let $X$ be a L\'evy process with no positive jumps and Laplace
exponent $f\in \B_3$. The L\'evy process reflected at its infimum,
$X_t-\inf_{0\leq s  \leq t} X_s$, is Markovian; and if $\tau$ denotes its
inverse local time at $0$, then  $\sigma=-X\circ \tau$ is a subordinator
called the  descending ladder-height
process. The Bernstein function
of the latter is then given by
 $g(q)=f(q)/q$; ; see e.g. Theorem VII.4(ii) in \cite{Belp}.

We next turn our attention to composition of functions;
here are some classical properties
\begin{proposition}\label{P3}
Consider two functions $f,g:]0,\infty[\to[0,\infty[$. Then we have the
implications
$$f,g\in\B_2 \ \Longrightarrow \ f\circ g\in \B_2\,,$$
$$f\in{\bf CM} \hbox{ and }g\in \B_2 \ \Longrightarrow \ f\circ g\in
{\bf CM}\,,$$
$$f\in\B_1 \hbox{ and }g\in \B_2 \ \Longrightarrow \ f\circ g\in
\B_1\,.$$
\end{proposition}
The first statement in Proposition \ref{P3} is related to the
celebrated subordination of Bochner (see, e.g. Section 3.9 in
\cite{Jacob} or Chapter 6 in
\cite{Sato}); more precisely if $\sigma$ and $\tau$ are two independent
subordinators with respective Bernstein functions $f_{\sigma}$ and
$f_{\tau}$, then $\sigma\circ \tau$ is again a subordinator whose
Bernstein function is $f_{\tau}\circ f_{\sigma}$. The second
statement is a classical result which can be found as Criterion 2 on page
441 in Feller \cite{Feller}; it is also related to Bochner's
subordination.

Finally we turn our attention to inverses.

\begin{proposition}\label{P4} Consider a
function $f:]0,\infty[\to]0,\infty[$. Then
$$f\in\B_2\cup \B_3 \ \Longrightarrow \ 1/f\in 
{\bf CM}\,.$$

Further, if $f^{-1}$ denotes the inverse of  $f$ when the latter is a
bijection, then
$$f\in\B_3\,,f\not\equiv 0 \ \Longrightarrow \ f^{-1}\in \B_2\,.$$

\end{proposition}

We mention that
if $f\in\B_3$, the completely monotone function $1/f$ is the Laplace
transform of the so-called scale function of the L\'evy process $X$ with
no positive jumps which has Laplace exponent $f$. See Theorem VII.8 in
\cite{Belp}. On the other hand,
$f^{-1}$ is the Bernstein function of the subordinator
of first-passage times $T_t:=\inf\left\{s\geq0: X_s>t\right\}$;
see e.g. Theorem VII.1 in \cite{Belp}. Finally, in the case when
$f\in\B_2$ is a Bernstein function, the completely monotone function $1/f$
is the Laplace transform of the renewal measure $U(dx)=\int_{0}^{\infty}
\P(\sigma_t\in dx)dt$, where $\sigma$ is a subordinator with Bernstein
function $f$.

\begin{corollary}\label{C2} Let $\Psi\not\equiv 0$ be a function in
$\B_3$, and denote by $\Phi=\Psi^{-1}\in\B_2$ its inverse bijection.
Then $q\to1/\Phi'(q)$ and ${\rm Id}/\Phi: q\to q/\Phi(q)$ are Bernstein
functions. Furthermore  ${1\over \Phi \Phi'}: q\to1/(\Phi(q) \Phi'(q))$ is
completely monotone.
\end{corollary}

\proof We know from Propositions
\ref{P1} and \ref{P4} that both $\Phi$ and $\Psi'$ are Bernstein
functions. We conclude from Proposition \ref{P3}
that $1/\Phi'=\Psi'\circ \Phi$ is again in $\B_2$.

Similarly, we know from Proposition \ref{P2} that
$q\to \Psi(q)/q$ is a Bernstein function, and  composition on the right
by the Bernstein function $\Phi$ yields ${\rm Id}/\Phi$ that is again 
in $\B_2$.

Finally, we can write $1/(\Phi \Phi')=f\circ \Phi$ where
$f(q)=\Psi'(q)/q$. We know from Proposition \ref{P1} that $\Psi'\in\B_2$,
so $f\in{\rm CM}$ by Proposition \ref{P2}. Since $\Phi\in\B_2$, we
conclude from Proposition \ref{P3} that $f\circ \Phi\in{\bf CM}$. \QED

If $\Phi=\Psi^{-1}$ is the Bernstein
function given by the inverse of a function $\Psi\in\B_3$, the Bernstein
function $1/\Phi'$ is the exponent of the subordinator $L^{-1}$ defined
as the inverse of the local time at $0$ of the L\'evy process with no
positive jumps and Laplace exponent $\Psi$. See e.g. Exercise VII.2 in
\cite{Belp}. On the other hand, ${\rm Id}/\Phi$ is then the Bernstein
function of the decreasing ladder times, see Theorem VII.4(ii) in
\cite{Belp}. The interested reader is also referred to \cite{Be1} for further
 factorizations for Bernstein functions which arise naturally for
L\'evy processes with no positive jumps, and their probabilistic
interpretations.

Next, recall that a function $f:]0,\infty[\to\R_+$ is called a Stieltjes
transform if it can be expressed in the form 
$$f(q)\,=\,b+\int_{[0,\infty[}{\nu(dt)\over t+q}\,,\qquad q>0\,,$$
where $b\geq0$ and  $\nu$ is a Radon measure on $\R_+$ such that 
$\int_{[0,\infty[}(1\wedge t^{-1})\nu(dt)<\infty$.
Equivalently, a Stieltjes
transform is the Laplace transform of a Radon
measure $\mu$ on $\R_+$ of the type $\mu(dx)=b\delta_0(dx)+h(x)dx$,
where $b\geq0$ and $h$ is a completely monotone function which belongs to
$L^1(\e^{-qx}dx)$ for every $q>0$; see e.g. Section 3.8 in \cite{Jacob}.

\begin{corollary}\label{C3} Let $f\in \B_2$ be a Bernstein
function such that its derivative $f'$ is a Stieltjes transform.
Then for every Bernstein function $g\in B_2$, the function
$f\circ {1\over g}$ is completely monotone.
\end{corollary}

\proof We can write
$$f(q)\,=\,a+ bq + \int_{0}^{q}dr\int_{0}^{\infty}dx\e^{-rx}h(x)\,,\qquad
q>0\,,$$
where $a,b\geq0$ and $h\in\B_1$. Thus 
$$f(q)\,=\,a+ bq + \int_{0}^{\infty}dx(1-\e^{-qx}){h(x)\over x}\,,\qquad
q>0\,,$$
and then
$$f\circ {1\over g}(q)\,=\,a+{b\over
g(q)}+\int_{0}^{\infty}dx(1-\e^{-x/g(q)}){h(x)\over x}\,.$$
We already know from Proposition \ref{P4} that $a+b/g \in {\bf CM}$.
The change of variable $y=x/g(q)$ yields
$$\int_{0}^{\infty}dx(1-\e^{-x/g(q)}){h(x)\over x}
\,=\,\int_{0}^{\infty}(1-\e^{-y})h(yg(q)){dy\over y}\,.$$
For each fixed $y>0$, $yg$ is a Bernstein function, so by Proposition
\ref{P3}, the function $q\to h(yg(q))$ is completely monotone.

We conclude that for every integer $n\geq
0$,
$$(-1)^n{\partial^n\over \partial q^n}(f\circ {1\over g})(q)
\,=\,\int_{0}^{\infty}(-1)^n{\partial^n\over \partial q^n}(h(yg(\cdot))(q)
(1-\e^{-y}){dy\over y}\,\geq\,0\,,$$
which establishes our claim. \QED

\end{document}